\documentclass[reqno]{amsart}
\usepackage{setspace}
\usepackage{graphicx}
\usepackage{amsmath}
\usepackage{amsfonts}
\usepackage{amssymb}

\begin{document}

\centerline{\LARGE \bf Twofold Translative Tiles in Three-Dimensional Space}

\bigskip\medskip
\centerline{\bf Mei Han, Qi Yang,\footnote{Both Qi Yang and Mei Han are the first author.} Kirati Sriamorn and Chuanming Zong}

\vspace{1cm}
\centerline{\begin{minipage}{12cm}
{\bf Abstract.} This paper proves the following statement: {\it If a convex body can form a twofold translative tiling in $\mathbb{E}^3$, it must be a parallelohedron.} In other words, it must be a parallelotope, a hexagonal prism, a rhombic dodecahedron, an elongated dodecahedron, or a truncated octahedron.
\end{minipage}}

\vspace{0.4cm}
\hspace{1.1cm}{2010 MSC: 52C22, 05B45, 52C17}

\vspace{1cm}
\noindent
{\Large\bf 1. Introduction}

\bigskip\noindent
Let $K$ be a convex body with interior $int(K)$ and boundary $\partial (K)$, and let $X$ be a discrete set, both in $\mathbb{E}^n$. We call $K+X$ a {\it translative tiling} of $\mathbb{E}^n$ and call $K$ a {\it translative tile} if $K+X=\mathbb{E}^n$ and the translates ${\rm int}(K)+{\bf x}_i$ are pairwise disjoint. In other words, if $K+X$ is both a packing in $\mathbb{E}^n$ and a covering of $\mathbb{E}^n$. In particular, we call $K+\Lambda$ a {\it lattice tiling} of $\mathbb{E}^n$ and call $K$ a {\it lattice tile} if $\Lambda $ is an $n$-dimensional lattice. Apparently, a translative tile must be a convex polytope. Usually, a lattice tile is called a {\it parallelohedron}.

In 1885, Fedorov \cite{fedo} characterized the two- and three-dimensional lattice tiles: {\it A two-dimensional lattice tile is either a parallelogram or a centrally symmetric hexagon; A three-dimensional lattice tile must be a parallelotope, a hexagonal prism, a rhombic dodecahedron, an elongated dodecahedron, or a truncated octahedron.}

The situations in higher dimensions turn out to be very complicated. Through the works of Delone \cite{delo}, $\check{S}$togrin \cite{stog} and Engel \cite{enge}, we know that there are exact $52$ combinatorially different types of parallelohedra in $\mathbb{E}^4$. A computer classification for the five-dimensional parallelohedra was announced by Dutour Sikiri$\acute{\rm c}$, Garber, Sch$\ddot{\rm u}$rmann and Waldmann \cite{dgsw} only in 2015.

Let $\Lambda $ be an $n$-dimensional lattice. The {\it Dirichlet-Voronoi cell} of $\Lambda $ is defined by
$$D=\left\{ {\bf x}: {\bf x}\in \mathbb{E}^n,\ \| {\bf x}, {\bf o}\|\le \| {\bf x}, \Lambda \|\right\},$$
where $\| X, Y\|$ denotes the Euclidean distance between $X$ and $Y$. Clearly, $D+\Lambda $ is a lattice tiling and the Dirichlet-Voronoi cell $D$ is a parallelohedron. In 1908, Voronoi \cite{voro} made a conjecture that {\it every parallelohedron is a linear image of the Dirichlet-Voronoi cell of a suitable lattice.} In $\mathbb{E}^2$, $\mathbb{E}^3$ and $\mathbb{E}^4$, this conjecture was confirmed by Delone \cite{delo} in 1929. In higher dimensions, it is still open.

To characterize the translative tiles is another fascinating problem. First it was shown by Minkowski \cite{mink} in 1897 that {\it every translative tile must be centrally symmetric}. In 1954, Venkov \cite{venk} and Aleksandrov \cite{alek} proved that {\it all translative tiles are parallelohedra.} Later, a new proof for this beautiful result was independently discovered by McMullen \cite{mcmu} (see also \cite{zong96}).

Let $X$ be a discrete multiset in $\mathbb{E}^n$ and let $k$ be a positive integer. We call $K+X$ a {\it $k$-fold translative tiling} of $\mathbb{E}^n$ and call $K$ a {\it $k$-fold translative tile} if every point ${\bf x}\in \mathbb{E}^n$ belongs to at least $k$ translates of $K$ in $K+X$ and every point ${\bf x}\in \mathbb{E}^n$ belongs to at most $k$ translates of ${\rm int}(K)$ in ${\rm int}(K)+X$. In other words, if $K+X$ is both a $k$-fold packing in $\mathbb{E}^n$ and a $k$-fold covering of $\mathbb{E}^n$. In particular, we call $K+\Lambda$ a {$k$-fold lattice tiling} of $\mathbb{E}^n$ and call $K$ a {\it $k$-fold lattice tile} if $\Lambda $ is an $n$-dimensional lattice. Apparently, a $k$-fold translative tile must be a convex polytope.

Multiple tilings were first investigated by Furtw\"angler \cite{furt} in 1936 as a generalization of Minkowski's conjecture on cube tilings. Let $C$ denote the $n$-dimensional unit cube. Furtw\"angler made a conjecture that {\it every $k$-fold lattice tiling $C+\Lambda$ has twin cubes. In other words, every multiple lattice tiling $C+\Lambda$ has two cubes which share a whole facet.} In the same paper, he proved the two- and three-dimensional cases. Unfortunately, when $n\ge 4$, this beautiful conjecture was disproved by Haj\'os \cite{hajo} in 1941. In 1979, Robinson \cite{robi} determined all the integer pairs $\{ n,k\}$ for which Furtw\"angler's conjecture is false. We refer to Zong \cite{zong05,zong06} for detailed accounts on this fascinating problem.

\medskip
Clearly, one of the most important and natural problems in multiple tilings is to classify or characterize the $n$-dimensional translative $k$-tiles and the $n$-dimensional lattice $k$-tiles (see Problems 1-4 at the end of Gravin, Robins and Shiryaev \cite{grs}). In the plane, it was proved by Yang and Zong \cite{yz1,yz2,zong20,zong-x} that, {\it besides parallelograms and centrally symmetric hexagons, there is no other two-, three- or fourfold translative tile. However, there are three classes of other fivefold translative tiles.} For more related results on multiple tilings, we refer to \cite{boll,gkrs,Grepstad,kolo,Lev-Liu}. As a counterpart of the above result and as a generalization of the theorem of Venkov, Aleksandrov and McMullen, this paper proves the following result.

\medskip\noindent
{\bf Theorem 1.} {\it Every three-dimensional twofold translative tile is a parallelohedron. In other words, every three-dimensional twofold translative tile must be a parallelotope, a hexagonal prism, a rhombic dodecahedron, an elongated dodecahedron, or a truncated octahedron.}

\medskip\noindent
{\bf Remark 1.} In 2016, K. Sriamorn \cite{sriamorn} announced a proof for this theorem in arXiv.org. Unfortunately, there was a gap in the last part of his proof.

\vspace{0.6cm}
\noindent
{\Large\bf 2. Preparations}

\medskip
\noindent
In this section, we recall some known concepts and results which will be useful for this paper.

\medskip
In 1885, E. S. Fedorov studied the two- and three-dimensional lattice tiles. He proved the following result.

\smallskip\noindent
{\bf Lemma 1 (Fedorov \cite{fedo}).} {\it A two-dimensional lattice tile is either a parallelogram or a centrally symmetric hexagon; A three-dimensional lattice tile must be a parallelotope, a hexagonal prism, a rhombic dodecahedron, an elongated dodecahedron, or a truncated octahedron.}

\smallskip
Tilings in higher dimensions have been studied by Minkowski \cite{mink}, Voronoi \cite{voro}, Delone \cite{delo}, Venkov \cite{venk}, Alexsandrov \cite{alek}, McMullen \cite{mcmu} and many others. Here we only introduce two basic results.

\smallskip
\noindent
{\bf Definition 1.} Let $P$ denote an $n$-dimensional centrally symmetric convex polytope with centrally symmetric facets and let $V$ denote a $(n-2)$-dimensional face of $P$. We call the collection of all those facets of $P$ which contain a translate of $V$ as a subface a belt of $P$.

\medskip\noindent
{\bf Lemma 2 (Venkov \cite{venk} and McMullen \cite{mcmu}).} {\it  A convex body $K$ is a translative tile if, and only if, it is a centrally symmetric polytope with centrally symmetric facets, such that each belt contains four or six facets.}

\medskip\noindent
{\bf Lemma 3 (Venkov \cite{venk} and McMullen \cite{mcmu}).}  {\it Every translative tile is a parallelohedron.}

\medskip\noindent
\smallskip
\noindent
{\bf Definition 2.} Let $P$ be an $n$-dimensional convex polytope. We call it a zonotope if it is a Minkowski sum of finite number of segments.
In other words,
$$P=S_1+S_2+\ldots +S_m,$$
where $S_i$ are segments in $\mathbb{E}^n$.

\medskip
In 2012, N. Gravin, S. Robins and D. Shiryaev studied multiple tilings in general dimensions and discovered the following result.

\medskip\noindent
{\bf Lemma 4 (Gravin, Robins and Shiryaev \cite{grs}).}  {\it An $n$-dimensional $k$-fold translative tile is a centrally symmetric polytope with centrally symmetric facets. In particular, in $\mathbb{E}^3$, every $k$-fold translative tile is a zonotope.}

\medskip
Recent years, Q. Yang and C. Zong systematically studied multiple translative tilings in the plane. One of their key ideas is the introduction of the adjacent wheel. This idea, by projection, is also useful in this paper.

Let $P_{2m}$ denote a centrally symmetric convex $2m$-gon centered at the origin, let ${\bf v}_1$, ${\bf v}_2$, $\ldots$, ${\bf v}_{2m}$ be the $2m$ vertices of $P_{2m}$ enumerated clock-wise, and let $G_1$, $G_2$, $\ldots $, $G_{2m}$ be the $2m$ edges, where $G_i$ is ended by ${\bf v}_i$ and ${\bf v}_{i+1}$. For convenience, we write
$$V=\{{\bf v}_1, {\bf v}_2, \ldots, {\bf v}_{2m}\}$$
and
$$\Gamma=\{G_1, G_2, \ldots, G_{2m}\}.$$

Assume that $P_{2m}+X$ is a $k$-fold translative tiling in $\mathbb{E}^2$, where $X=\{{\bf x}_1, {\bf x}_2, {\bf x}_3, \ldots \}$ is a discrete multiset with ${\bf x}_1={\bf o}$. Now, let us observe the local structures of $P_{2m}+X$ at the vertices ${\bf v}\in V+X$.

Let $X^{\bf v}$ denote the subset of $X$ consisting of all points ${\bf x}_i$ such that
$${\bf v}\in \partial (P_{2m})+{\bf x}_i.$$
Since $P_{2m}+X$ is a multiple tiling, the set $X^{\bf v}$ can be divided into disjoint subsets $X^{\bf v}_1$, $X^{\bf v}_2$, $\ldots ,$ $X^{\bf v}_t$ such that the translates in $P_{2m}+X^{\bf v}_j$ can be re-enumerated as $P_{2m}+{\bf x}^j_1$, $P_{2m}+{\bf x}^j_2$, $\ldots $, $P_{2m}+{\bf x}^j_{s_j}$ satisfying the following conditions:

\medskip
\noindent
{\bf 1.} {\it ${\bf v}\in \partial (P_{2m})+{\bf x}^j_i$ holds for all $i=1, 2, \ldots, s_j.$}

\smallskip\noindent
{\bf 2.} {\it Let $\angle^j_i$ denote the inner angle of $P_{2m}+{\bf x}^j_i$ at ${\bf v}$ with two half-line edges $L^j_{i,1}$ and $L^j_{i,2}$, where $L^j_{i,1}$, ${\bf x}^j_i-{\bf v}$ and $L^j_{i,2}$ are in clock order. Then, the inner angles join properly as
$$L^j_{i,2}=L^j_{i+1,1}$$
holds for all $i=1,$ $2,$ $\ldots ,$ $s_j$, where $L^j_{s_j+1,1}=L^j_{1,1}$.}

\medskip
For convenience, we call such a sequence $P_{2m}+{\bf x}^j_1$, $P_{2m}+{\bf x}^j_2$, $\ldots $, $P_{2m}+{\bf x}^j_{s_j}$ an {\it adjacent wheel} at ${\bf v}$. In other words, if ${\bf v}$ belongs to the boundary of a tile then we follow this tile around, moving from tile to tile, until it closes up again. It is easy to see that
$$\sum_{i=1}^{s_j}\angle^j_i =2w_j\cdot \pi$$
hold for positive integers $w_j$. Then we define
$$\varpi ({\bf v})=\sum_{j=1}^tw_j= {1\over {2\pi }}\sum_{j=1}^t\sum_{i=1}^{s_j}\angle^j_i\eqno(1)$$
and
$$\varphi ({\bf v})=\sharp \left\{ {\bf x}_i:\ {\bf x}_i\in X,\ {\bf v}\in {\rm int}(P_{2m})+{\bf x}_i\right\}.\eqno(2)$$
In other words, $\varpi ({\bf v})$ is the tiling multiplicity produced by the boundary and $\varphi ({\bf v})$ is the tiling multiplicity produced by the interior.

\medskip
Clearly, if $P_{2m}+X$ is a $k$-fold translative tiling of $\mathbb{E}^2$, then
$$k= \varphi ({\bf v})+\varpi ({\bf v})\eqno (3)$$
holds for all ${\bf v}\in V+X$.

\medskip\noindent
{\bf Lemma 5 (Yang and Zong \cite{yz2}).} {\it Assume that $P_{2m}$ is a centrally symmetric convex $2m$-gon centered at the origin, $P_{2m}+X$ is a translative multiple tiling of the plane, and $\mathbf{v}\in V+X$. Then we have
$$\varpi ({\bf v})=\kappa\cdot {{m-1}\over 2}+\ell\cdot {1\over 2},$$
where $\kappa $ is a positive integer and $\ell$ is the number of the edges in $\Gamma+X$ which take ${\bf v}$ as an interior point.}

\vspace{0.6cm}
\noindent
{\Large\bf 3. A Proof for Theorem 1}

\medskip
Now we restrict ourselves in the three-dimensional Euclidean space.

\medskip\noindent
{\bf Lemma 6.} {\it Suppose that $P$ is a twofold translative tile in $\mathbb{E}^{3}$. Let $G$ be an edge (one-dimensional face) of $P$ and let $B(G)$ be the belt determined by $G$, which consists of $2m$ facets $F_{1}$, $F_2$, $...,$ $F_{2m}$. Let $G_{1},$ $G_2$, $...,$ $G_{2m}$ be the translates of $G$ such that $G_{i}, G_{i+1}\subset F_{i}$ and $G=G_{1}$. Suppose that $G_{i+1}=G_{i}+\mathbf{g}_{i}$, where $\mathbf{g}_{i}\in\mathbb{E}^{3}$ and $1<i<m-1$, then $rint(G_{1})+\mathbf{g}_{i}\subset int(P)$.}

\medskip\noindent{\bf Proof.} According to Lemma 4, the twofold translative tile $P$ is a centrally symmetric polytope with centrally symmetric facets. In fact, it is a zonotope. Then, without loss of generality, we have
$$P=S_1+S_2+\ldots +S_m,$$
where $S_i$ are segments centered at the origin of $\mathbb{E}^{3}$ and no pair of them are colinear.

Without loss of generality, we assume that $G_1$ is a translate of $S_1$, which is in the direction of ${\bf e}_3=(0,0,1)$, and $F_i$ is a translate
of
$$Q_i=S_1+S_2+\ldots +S_n,\eqno(4)$$
where $n<m$. Clearly, all $S_1$, $S_2$, $\ldots ,$ $S_n$ are coplanar and no other segment coplanar with them. Then, we write
$$Q'_i=S_{n+1}+S_{n+2}+\ldots +S_m.\eqno(5)$$
Clearly, we have
$$P=Q_i+Q'_i\eqno(6)$$
and
$$F_i=Q_i+{\bf c}_i,$$
where ${\bf c}_i\in Q'_i$. Furthermore, if
$$G_1\subset Q_i+{\bf d}_i,\qquad {\bf d}_i\in Q'_i,$$
then we have
$$Q_i+{\bf d}_i\subset P$$
and therefore
$$G_1+{\bf g}_i\subset Q_i+{\bf d}_i\subset P.\eqno(7)$$

Assume that ${\bf x}_i$ is a point such that
$$G_1=G_{i+m+1}+{\bf x}_i.$$
Then we have
$$-F_i+{\bf x}_i=Q_i+{\bf d}_i.\eqno(8)$$
Since $1<i<m-1$, by projecting to $xy$-plane (see Figure 1), it is easy to see that
$$ int(P\cap (P+{\bf x}_i))\not= \emptyset.\eqno(9)$$

\begin{figure}[ht]
\centering
\includegraphics[height=5.5cm,width=8.8cm,angle=0]{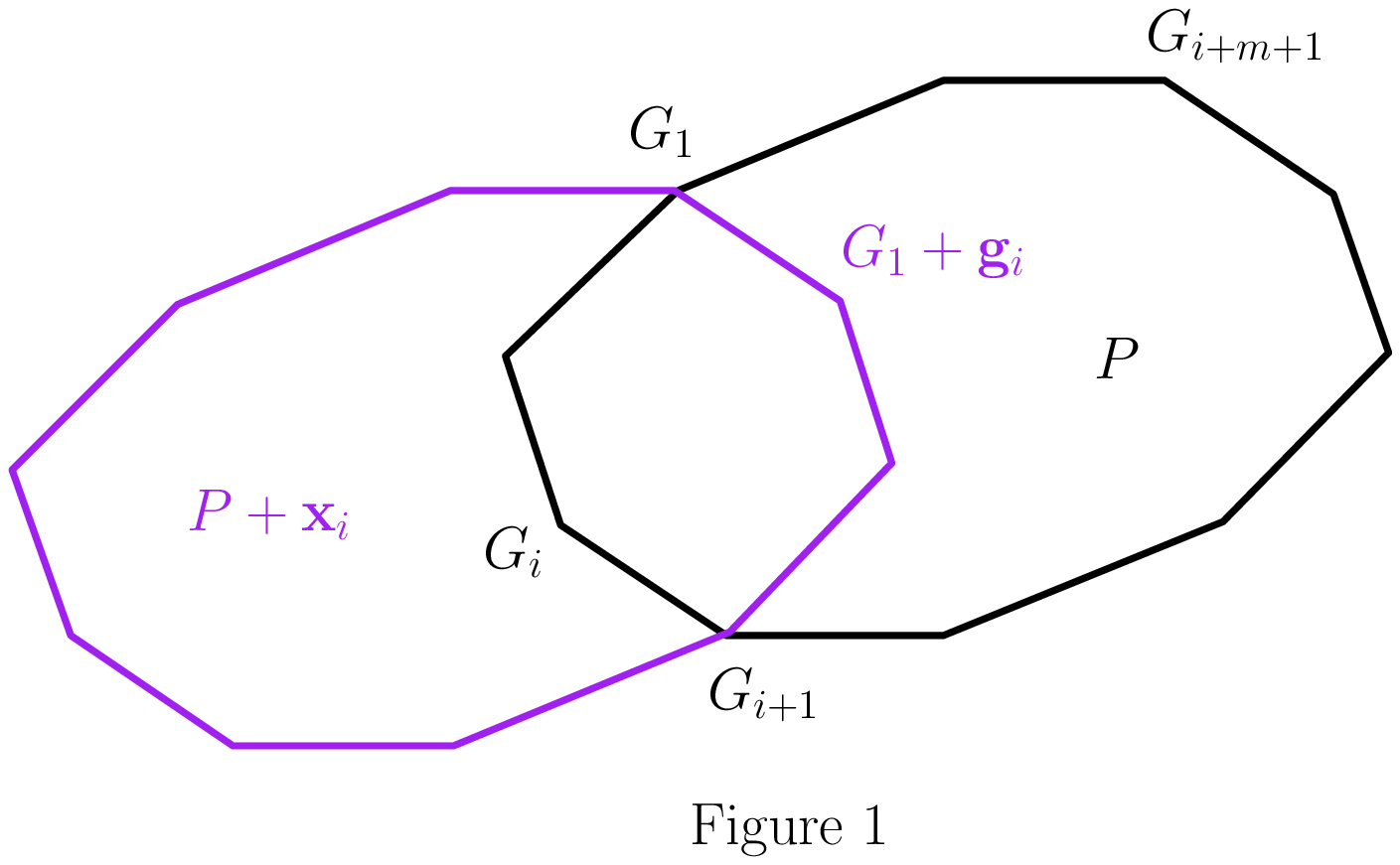}
\end{figure}

\noindent
On the other hand, it is known that $P\cap (P+{\bf x}_i)$ is centrally symmetric. If
$$rint (G_1)+{\bf g}_i\not\subset int(P),\eqno(10)$$
it can be deduced from (7) that the whole segment $G_1+{\bf g}_i$ must belong to $\partial (P)$. When we project $P\cap (P+{\bf x}_i)$ onto the $xy$-plane, we obtain a centrally symmetric polygon, half from the boundary of $P$ and the other half from $P+{\bf x}_i$. But now the half from $P+{\bf x}_i$ is only one segment. Thus, we have
$$ int(P\cap (P+{\bf x}_i)) = \emptyset,\eqno(11)$$
which contradicts to (9). As a conclusion, we must have
$$rint (G_1)+{\bf g}_i\subset int(P).\eqno(12)$$
The Lemma is proved. \hfill{$\Box$}

\medskip
\noindent
{\bf Proof of Theorem 1.} Assume that $P$ is a twofold translative tile and $P+X$ is a corresponding twofold tiling in $\mathbb{E}^{3}$. Assume that $G$ is an edge of $P$ such that $B(G)$ has $2m$ facets $F_1$, $F_2$, $\ldots,$ $F_{2m}$.

Similar to (1) and (2), for any point ${\bf v}\in rint(G)$, we define $\varpi ({\bf v})$ and $\varphi ({\bf v})$. Of course, here $\varpi ({\bf v})$ is for the corresponding dihedral angle. Since the twofold translative tiling is invariant under linear transformations, we may assume further that $G$ is perpendicular to the $xy$-plane. Now, we consider three cases.

\smallskip\noindent
{\bf Case 1.} $m\ge 6$. Then, it follows by Lemma 5 and (3) that
$$\varpi ({\bf v})=\kappa\cdot {{m-1}\over 2}+\ell\cdot {1\over 2}>2\eqno(13)$$
and
$$\varphi ({\bf v})+\varpi ({\bf v})>2,\eqno(14)$$
which contradicts the assumption that $P+X$ is a twofold tiling in $\mathbb{E}^{3}$.

\smallskip\noindent
{\bf Case 2.} $m=5$. Then, it follows from Lemma 5 that the only possibility is both $\varpi ({\bf v})=2$ and $\varphi ({\bf v})=0$ hold for all ${\bf v}\in G$. Then, by projecting all the translates around ${\bf v}$ in $P+X$ onto the $xy$-plane we obtain the following configuration, where $P_{10}$ is the projection of $P$.

\begin{figure}[ht]
\centering
\includegraphics[height=6.6cm,width=10cm,angle=0]{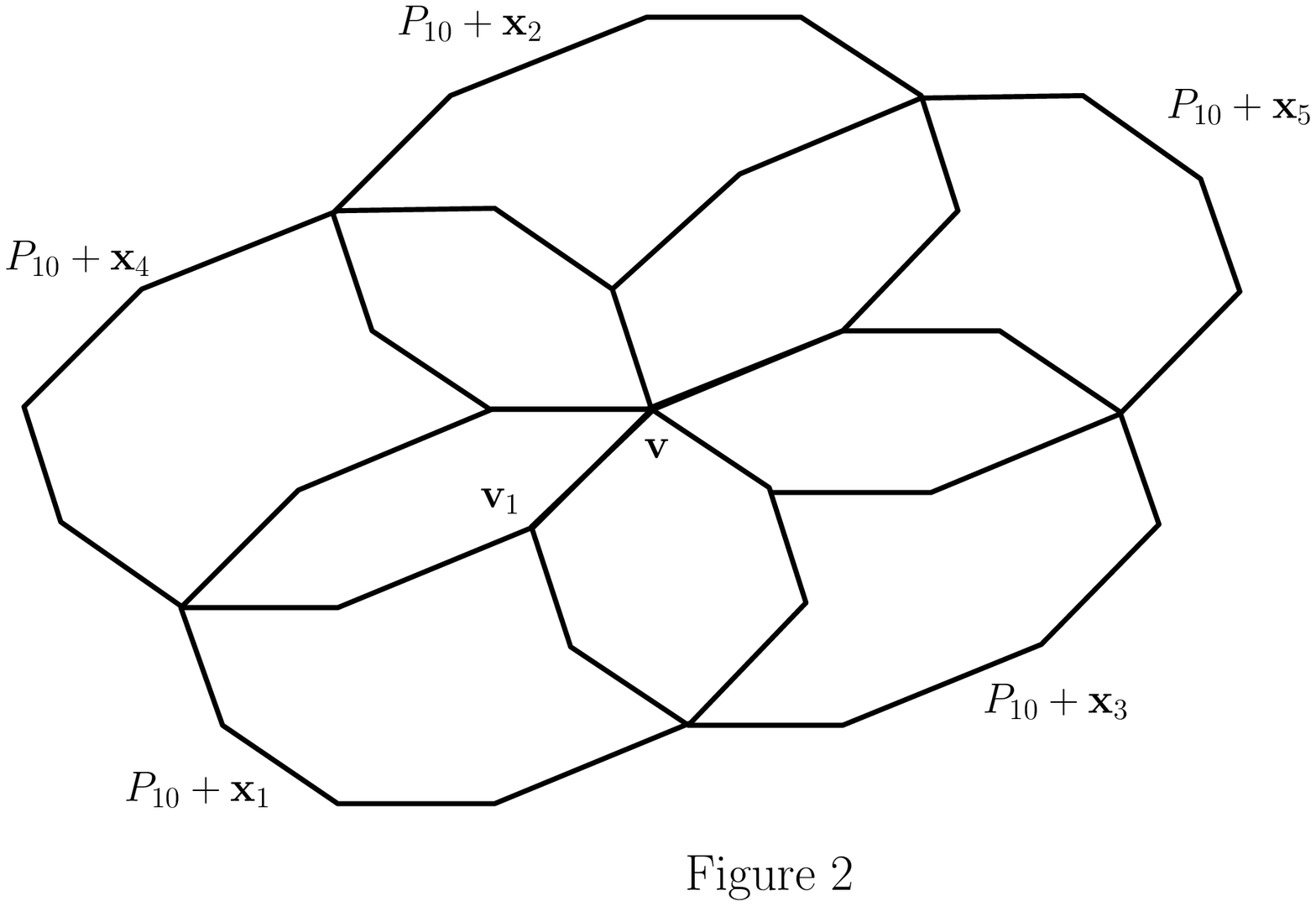}
\end{figure}

Assume that ${\bf v}_1$ is the projection of the vertical edge $G'$ of $P+{\bf x}_3$, which is a translate of $G$. It follows from Lemma 6 that
$$rint(G')\cap (int(P)+{\bf x}_1)\not= \emptyset\eqno(15)$$
and therefore
$$\varphi ({\bf v}')\ge 1\eqno(16)$$
holds for some ${\bf v}'\in G'$. Then, applying Lemma 5 to ${\bf v}'$, we obtain
$$\varpi ({\bf v}')\ge 2.\eqno(17)$$
Consequently, by (3), (16) and (17) we get
$$\varphi ({\bf v}')+\varpi ({\bf v}')\ge 3,$$
which contradicts the assumption that $P+X$ is a twofold tiling in $\mathbb{E}^{3}$.

\smallskip\noindent
{\bf Case 3.} $m=4$. Then, it follows from Lemma 5 that the only possibility is both $\varpi ({\bf v})=2$ and $\varphi ({\bf v})=0$ hold for all ${\bf v}\in G$. Then, by projecting all the translates around ${\bf v}$ in $P+X$ onto the $xy$-plane we obtain the following configuration, where $P_8$ is the projection of $P$.

\begin{figure}[ht]
\centering
\includegraphics[height=7cm,width=8cm,angle=0]{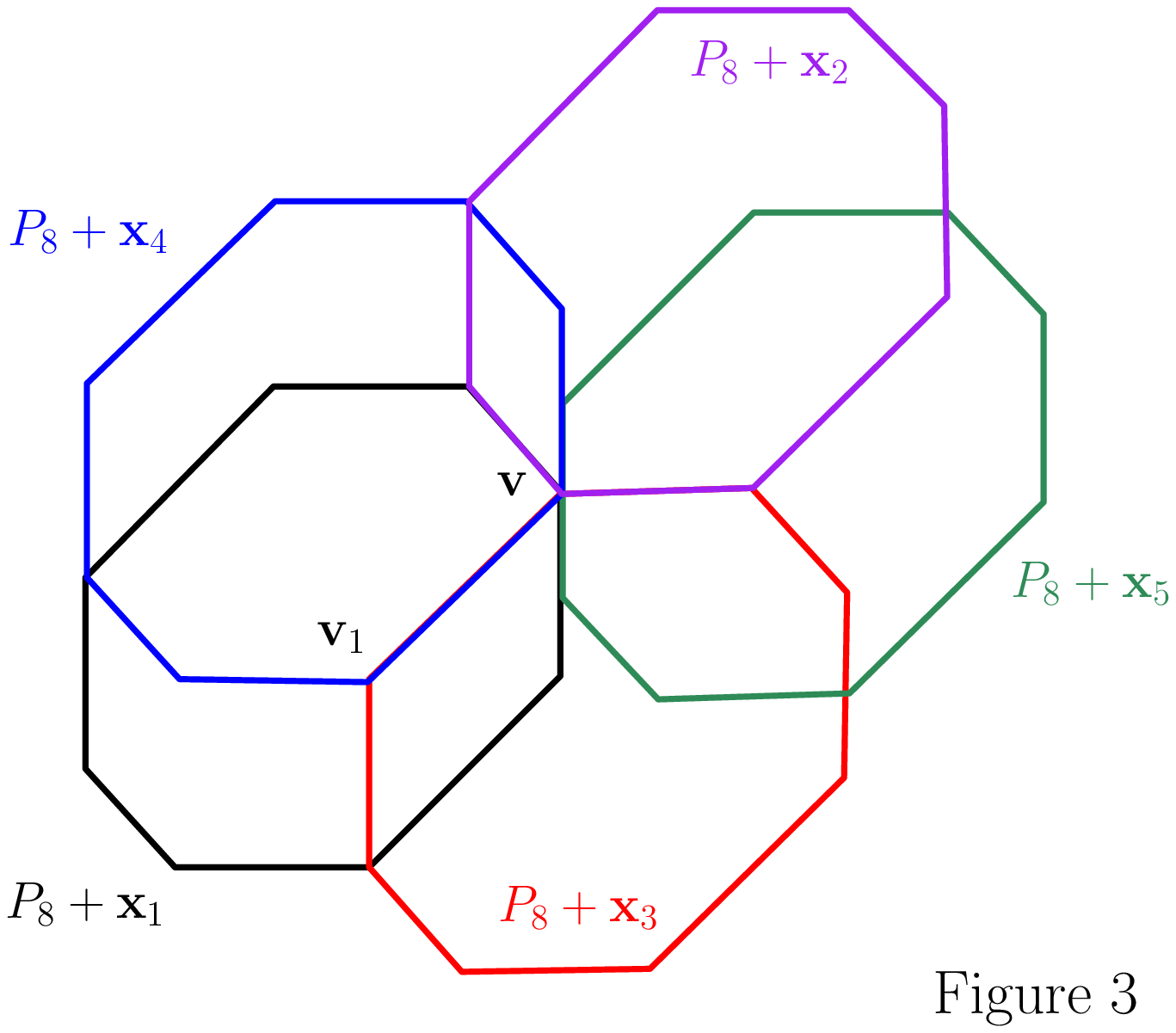}
\end{figure}

Assume that ${\bf v}_1$ is the projection of the vertical edge $G'$ of $P+{\bf x}_3$, which is a translate of $G$. It follows from Lemma 6 that
$$rint(G')\cap (int(P)+{\bf x}_1)\not= \emptyset\eqno(19)$$
and therefore
$$\varphi ({\bf v}')\ge 1\eqno(20)$$
holds for some ${\bf v}'\in G'$. Then, applying Lemma 5 to ${\bf v}'$, we obtain
$$\varpi ({\bf v}')\ge {3\over 2}.\eqno(21)$$
Consequently, by (3), (20) and (21) we get
$$\varphi ({\bf v}')+\varpi ({\bf v}')\ge 3,\eqno(22)$$
which contradicts the assumption that $P+X$ is a twofold tiling in $\mathbb{E}^{3}$.

\medskip
As a conclusion of these three cases, we get $m\le 3$ for every belt of $P$. Then, the theorem follows easily from Lemma 2, Lemma 3 and Lemma 1. \hfill{$\Box$}

\vspace{0.6cm}\noindent
{\bf Acknowledgements.} This work is supported by ERC Starting Grant No. 713927, the National Natural Science Foundation of China (NSFC11921001), the National Key Research and Development Program of China (2018YFA0704701), and 973 Program 2013CB834201.

\vspace{0.6cm}
\noindent
Mei Han, Center for Applied Mathematics, Tianjin University, Tianjin 300072, China.

\noindent
Qi Yang, Department of Mathematics, Bar Ilan University, Israel.

\noindent
Kirati Sriamorn, Department of Mathematics and Computer Science, Chulalongkorn University, Thailand.

\noindent
Chuanming Zong, Center for Applied Mathematics, Tianjin University, Tianjin 300072, China.

\noindent
Email: cmzong@math.pku.edu.cn

\end{document}